\documentclass[a4paper]{amsart}
\usepackage{amssymb} 
\usepackage{graphicx}

\usepackage{color}
\newtheorem{theo}{Theorem}[section]
\newtheorem{prop}{Proposition}[section]

\newtheorem{lemma}{Lemma}[section]

\numberwithin{equation}{section} 

\def\t{\textrm}

\def\beq{\begin{eqnarray}}
\def\eeq{\end{eqnarray}}
\def\baa{\begin{array}}
\def\eaa{\end{array}}

\newcommand{\bdef}{\begin{definition}}
\newcommand{\be}{\begin{equation}}
\newcommand{\ee}{\end{equation}}
\newcommand{\bt}{\begin{theo}}
\newcommand{\et}{\end{theo}}
\newcommand{\bl}{\begin{lemma}}
\newcommand{\el}{\end{lemma}}
\newcommand{\bp}{\begin{prop}}
\newcommand{\epr}{\end{prop}}

\def\dsp{\displaystyle}

\def\ep{\epsilon}
\def\ov{\overline}

\let\dsp=\displaystyle
\def\R{{\mathbb R}}

\def\li{\lambda_i}
\def\lai{\lambda_i}
\def\v{v_i}
\def\u{u_i}
\def\vx{{v_i}_x}
\def\ux{{u_i}_x}

\def\vt{{v_i}_t}
\def\ut{{u_i}_t}

\def\beti{\beta_i}
\def\Di{D_i}
\def\phii{\phi_i}
\def\phx{{\phi_i}_x}
\def\pht{{\phi_i}_t}
\def\phxx{{\phi_i}_{xx}}
\def\phtx{{\phi_i}_{tx}}

\def\pho{\phi_{0i}}

\def\uo{u_{0i}}
\def\vo{v_{0i}}

\def\hu{\hat u}
\def\hv{\hat v}
\def\hf{\hat \phi}
\def\M{{\mathcal{M}}}
\def\U{\mathcal{O}}
\def\E{\mathcal{I}}
\def\A{\mathcal{A}}
\def\Ii{I_i}

\def\iIi{\int_{\Ii}}
\def\nld{\Vert_{2}}

\def\nhu{\Vert_{H^1}}

\def\soT{\sup_{[0,T]}}

\def\nhd{\Vert_{H^2}}
\def\ioT{\int_0^T}

\def\suM{\sum_{i\in \M}}

\def\suEn{\sum_{i\in \E^\nu}}
\def\suUn{\sum_{i\in \U^\nu}}
\def\suMj{\sum_{j\in \M}}
\def\suMjn{\sum_{j\in \M^\nu}}

\def\iota{\int_0^\tau}

\def\oespxi{\exp\left(\frac {{\phii}^0(x)} \li\right)}
\def\oespxj{\exp\left(\frac {{\phi}^0_j(x)} {\lambda_j}\right)}

\def\espxi{\exp\left(\frac {{\phii}(x)} \li\right)}

\def\espnju{\exp\left(\frac {{\phii}(N_\nu)} \li\right)}
\def\espniu{\exp\left(\frac {{\phi}_j(N_\nu)} {\lambda_j}\right)}

\def\ooespniu{\exp\left(\frac {\phi^0_j(N_\nu)} {\lambda_j}\right)}
\def\ooespnju{\exp\left(\frac {{\phii}^0(N_\nu)} \li\right)}

\begin{document}
\title[Stationary solutions and asymptotic behaviour]
{Stationary solutions  and asymptotic behaviour for a chemotaxis hyperbolic model on a network}
\author[ F. R. Guarguaglini]{ F. R. 
Guarguaglini$^\diamond$}
 
\thanks{ \noindent 
$\diamond$ Disim, Universit\'a degli
Studi di L'Aquila, Via Vetoio, I--67100 Coppito (L'Aquila), Italy. E--mail:
guarguag@univaq.it }
\date{}
\begin{abstract} 
This paper approaches the  question of  existence and uniqueness of stationary solutions  to
a semilinear hyperbolic-parabolic system  and the study of the asymptotic behaviour of  global  solutions.  
The system  is a model for  some biological
phenomena evolving on a network composed by  a finite number of nodes and  oriented arcs. 
The transmission conditions for the unknowns, set at each inner node, 
 are  crucial features  of the model.

\end{abstract}
\maketitle
\bigskip

\begin{section}
{\bf Introduction}

In this paper we consider a semilinear  hyperbolic-parabolic system  evolving on a finite planar network composed from  nodes connected by $m$ oriented arcs $I_i$, in one space dimension,

\begin{equation}\label{is}
    \left\{ \begin{array}{l}
\partial_t \u +\lambda_i \partial_x \v=0\ ,
\\ \\
\partial_t \v +\lambda_i \partial_x \u= \u \partial_x \phi_i -\beta_i \v\ ,
\qquad\qquad     t\geq 0\, , x\in I_i\, ,\  i=1,...,m;
\\ \\
\partial_t \phii= D_i \partial_{xx} \phi_i +a_i\u -b_i\phi_i\ ,
\end{array}\right.\end{equation}
the system is complemented by initial, boundary and transmission conditions at the nodes (see Section 2).

We are interested in the study of stationary solutions  and  asymptotic behaviour of global solutions of the problem.

The above system    has been proposed as a model for chemosensitive movements of bacteria or cells on an artificial scaffold \cite{noi}. The unknown $u$ stands for the cells concentration, $v$ is the average flux and $\phi$ is the chemo-attractant concentration. 
In particular, the model turns out to be useful to describe the process of dermal wound healing, when the stem cells in charge of the reparation of dermal tissue  (fibroblasts) create an extracellular matrix and move along it  to fill the wound, driven by chemotaxis;
 tissue engineers use artificial scaffolds, constituted by a network of crossed polymeric threads,  inserting them within the wound to accelerate the process (see
\cite{harley,mandal,spadac}).
 In the above mathematical model, the arcs of the graph  mimic the fibers of the scaffold; each of them is characterized by a tipical velocity $\li$, a friction coefficient $\beta_i$, a diffusion coefficient $D_i$, and a production rate $a_i$ and  a degradation one  $b_i$ ; the functions $\u,\phii$ are the densities of fibroblasts and chemoattractant  on each arc.

 Starting from the Keller-Segel paper \cite{KS} in 1970 until now, a lot of articles have been devoted  to PDE models in domains of $\R^n$ for chemotaxis phenomena. The parabolic (or parabolic-elliptic) Patlak-Keller-Segel system
 is the most studied model \cite{horst,perth,mur}; in recent years, hyperbolic models have been introduced too, in order to avoid the unrealistic infinite speed of propagation of cells, occurring in parabolic models \cite{DH,FLP,Hillen, perth,Segel,AG,hillerho,hillenst,gumanari}.
 
In \cite{gumanari} the Cauchy  and the Neumann problems  for the system in (\ref{is}), respectively in $\R$ and in bounded intervals of $\R$, are studied, providing existence of global solutions and stability of constant states results.

Recently an interest in these mathematical models evolving on networks is  arising, due to their applications in the study  of biological phenomena and traffic flows, both in  parabolic cases \cite{4,CC, mugn} and in hyperbolic ones \cite{pic,zua1,zua2,noi,BNR}.

We notice that the transmission conditions for the unknowns, at each inner node, which complement the equations on networks, are  crucial characteristics  of the model, since they are the coupling among the solution's components on each arc.

Most of the studies carried out until now, consider  continuity  conditions at each inner node for the density functions \cite{zua1,CC, mugn}; nevertheless, the eventuality of discontinuities  at the nodes seems   a  more appropriate framework to decribe 
 movements of individuals or traffic flows phenomena \cite{cgp}. 

For these reason in  \cite{noi},  transmission conditions which link the values of the density functions at the nodes with the fluxes, without imposing any continuity, are introduced; these conditions guarantee the fluxes conservation at each inner node, and,  at the same time, the m-dissipativity of the linear spatial differential operators, a crucial property in the proofs of existence of local and global solutions contained in that paper.

In this paper we focus our attention  on the  stationary solutions to problem (\ref{is}) complemented by null fluxes boundary conditions and by the same transmission conditions of 
\cite{noi} (see next section and Section 3 in \cite{noi} for details). 

In the case of acyclic networks, although the transmission conditions do not set the continuity of the density $u$  at the inner nodes, the fluxes conservation at those nodes  and the boundary null fluxes conditions  imply the absence of jumps discontinuities at the inner vertexes, for the component $u$ of a stationary solution.

On the other hand,  if the ratio between $a_i$ and $b_i$ is constant, it is easy to prove that, for any fixed  mass $\dsp \suM\iIi \u(x) dx$, a stationary solution  constant on  the whole network  exists and the constant values of the densities are  determined  by the mass.

Such constant states turn out to be  asymptotic profiles for a class of solutions; so, we point out that, for small global solutions to our problem,
the discontinuities at the inner nodes vanish when $t$ goes to infinity, since their asymptotic profiles are continuous functions on the whole network. 

In the next section we give the statement of the problem; in particular we introduce  the transmission conditions .

In Section 3 we 
consider acyclic networks  and we prove the existence and uniqueness of the  stationary solution when the  mass  is suitably small; if the parameters $a_i, b_i$  have constant ratio, such solution is a constant state, continuous on the network.
Moreover, for general networks,  we prove that the density $u$ of a stationary solution which is small in a suitable norm, has to be continuous.

The last section is devoted to study the asymptotic behaviour of solutions corresponding to initial data which are  small perturbations of small constant states. The results obtained in this section constitute  the sequel and the development of the result of existence of global solutions in \cite{noi} and the proofs are based on the same thecniques; in particular, 
simple supplements to the  a priori estimates for the solutions to (\ref{is}) obtained in \cite{noi}, allow to prove that the constant states are the asymptotic profiles.

The  study of the asymptotic behaviour provide informations about the evolution of  a small mass of individuals moving on a network driven by chemotaxis:  suitable   
initial distributions of individuals and chemoattractant, for large time  evolve towards constant  distributions  on the network, preserving the mass of individuals.

We recall that the stability of the constant solutions to this system, considered on bounded interval in $\R$, is studied in \cite{gumanari} and stationary solutions and asymptotic behaviour for a linear system of uncoupled conservation laws on network are studied in \cite{KZ}.
 
 Finally, in \cite{BNR} the authors introduce  a numerical scheme to approximate the solutions to the  problem  (\ref{sysi}); in that paper transmission conditions are set for  the Riemann invariants of the hyperbolic part of the system, $w_i^{\pm}=\frac 1 2 (u_i\pm v_i)$, 
 and are equivalent to our ones  for some choices of the transmission coefficients. The tests presented there, in the case of acyclic graph  and  dissipative transmission coefficients, show an asymptotic behaviour of the solutions which agrees with our theoretical results.

\end{section}

\begin{section} 
{\bf Statement of the evolution problem } 

We consider a finite connected graph $\mathcal G=(\mathcal V, \mathcal A)$ composed by a 
set  $\mathcal V$ 
of $n$ nodes (or vertexes)
 and a set $\mathcal A$
of $m$ oriented arcs, $\mathcal A=\{I_i:i\in \mathcal M=\{1,2,...,m\}\}$. 

Each node 
 is a point of the plane and 
each oriented arc $I_i$ is an oriented segment  joining two nodes.

We use   $e_j$,  $j\in \mathcal J$,  to indicate the external vertexes (or boundary vertexes) of the graph, i.e. the vertexes belonging  to only one  arc, and by
$I_{i(j)}$
  the  external arc  outgoing or incoming  in the external vertex $e_j$. 
  
 Moreover,  we use $N_\nu$, $\nu\in\mathcal P$, to denote the inner nodes; for each of them we consider the set of incoming arcs 
$\mathcal A_{in}^\nu=\{I_i:i\in \mathcal I^\nu\}$
and the set of the outgoing ones
$\mathcal A_{out}^\nu=\{I_i:i\in \mathcal O^\nu\}$; let $\mathcal M^\nu=\mathcal I^\nu\cup\mathcal O^\nu$.

In this paper, a  {\it path} in the  graph is  a sequence of arcs, two by two adjacent,
without taking into account  orientations . Moreover, we call {\it acyclic} a graph which does  not contain cycles:  for each couple of nodes, there exists a unique  path with no repeated arcs connecting them (an example of acyclic graph is in Fig. 1).

Each arc $I_i$ is considered as a compact one dimensional interval $[0,L_i]$.
A function $f$ defined on $\mathcal A$ is a m-tuple of functions $f_i$, $i\in \mathcal M$, each one defined 
on $I_i$;  $ f_i(N_\nu)$  denotes $f_i(0) $ if $N_\nu$ is the initial point of the arc $I_i$ and $f_i(L_i)$ if $N_\nu$ is the end point, and similarly for $f(e_j)$.

We set $\dsp L^p(\mathcal A):=\{f: f_i\in L^p(I_i), i\in \M\}$, 
$\dsp H^s(\mathcal A):=\{f :f_i\in H^s(I_i), i\in\M\}$
 and
$$\Vert f\Vert_2:=\suM \Vert f_i\Vert_2\ , \ 
\Vert f\Vert_{H^s}:=\suM \Vert f_i\Vert_{H^s}.$$ 

We consider the evolution  of the following one-dimensional problem on the graph $\mathcal G$
\begin{equation}\label{sysi}
    \left\{ \begin{array}{l}
\partial_t \u +\lambda_i \partial_x \v=0\ ,
\\ \\
\partial_t \v +\lambda_i \partial_x \u= \u \partial_x \phi_i -\beta_i \v\ ,
\qquad     t\geq 0  , \ x\in I_i  , \  i\in \M,
\\ \\
\partial_t \phii= D_i \partial_{xx} \phi_i +a_i\u -b_i\phi_i\ ,
\end{array}\right.\end{equation}
where $ a_i\geq 0\ ,\li\ b_i,D_i,\beta_i >0 $  .

We complement the system with the initial conditions
\be\label{uvic}
u_{i0},v_{i0}\in H^1(I_i)\ ,\ 
\phi_{i0}\in H^2(I_i)
\ \t{ for } i\in\M\ ;\ee

the boundary conditions  at each outer point $e_j$ are the null flux conditions

\be\label{vbc}
v_{i(j)}(e_j, t)=0\ ,
\qquad 
\  t>0\ , \ j\in\mathcal J\  ,
\ee
\be\label{phibc}
{\phi_{i(j)}}_x(e_j, t)=0\qquad 
\ t>0 \ ,\  j\in \mathcal J .\ee



In addition, at each inner node $N_\nu$ we impose the following transmission conditions for the unknown $\phi$
\be\label{KC}\left\{\begin{array}{ll}\displaystyle
\Di \phx(N_\nu,t)=  \sum_{j\in\M^\nu} \alpha_{ij}^\nu(\phi_j(N_\nu,t)-\phii(N_\nu,t))\  \ ,\ i\in \E^\nu \ ,\  t> 0\ ,
\\ \\ \dsp
-\Di \phx(N_\nu,t)=  \sum_{j\in\M^\nu} \alpha_{ij}^\nu(\phi_j(N_\nu,t)-\phii(N_\nu,t))\  ,\ \  i\in \U^\nu \ ,\  t> 0\ ,
\\ \\
\alpha_{ij}^\nu\geq 0\ ,\ \alpha_{ij}^\nu=\alpha_{ji}^\nu \ \t{ for all } i,j\in \M^\nu\ ,
\end{array}\right.\ee
which imply the continuity of the flux  at  each node, for all $t>0$,
$$\suEn D_i\phx(N_\nu,t)=\suUn D_i \phx(N_\nu,t)\ .$$

For the unknonws $v$ and $u$ we impose the transmission conditions
\be\label{TC}\left\{\begin{array}{ll}\dsp
-\lai\v(N_\nu,t)=\sum_{j\in\M^\nu}
K_{ij}^\nu \left(u_j(N_\nu,t)-\u(N_\nu,t)\right)\ ,\ i\in \E^\nu\ ,\ t>0\ ,\\ \\ \dsp
\lai\v(N_\nu,t)=\sum_{j\in\M^\nu}
K_{ij}^\nu \left(u_j(N_\nu,t)-\u(N_\nu,t)\right)\ ,\  i\in\U^\nu \ , \ t>0\ ,\\ \\
K_{ij}^\nu\geq 0\ ,\ K_{ij}^\nu=K_{ji}^\nu \ \t{ for all } i,j\in \M^\nu\ .
\end{array}\right.\ 
\ee
These conditions ensure the conservation of the 
flux of the density of cells at each node $N_\nu$ , for $t>0$,

$$
\suEn \lai\v(N_\nu,t)=\suUn\lai\v(N_\nu,t)\ ,
$$
which corresponds to the conservation of the total mass
$$\suM\iIi\u(x,t)\ dx=\suM\iIi\uo(x)\ dx\ ,$$
i.e. no death nor birth of individuals occours during the observation.

Motivations for the constraints on the coefficients in the transmission conditions can be found in \cite{noi} .

Finally, we impose the following compatibility conditions
 \be\label{compat} 
u_{i0},v_{i0}, \phi_{i0}\t{  satisfy conditions (\ref{vbc})-(\ref{TC})  
for all } i\in\M\ .
\ee 

Existence and uniqueness of  local solutions 
 to problem (\ref{sysi})-(\ref{compat}),
$$u,v \in C([0,T];H^1(\mathcal A))\cap C^1([0,T];L^2(\mathcal A))\ ,\ 
\phi \in C([0,T];H^2(\mathcal A))\cap C^1([0,T];L^2(\mathcal A))\ $$ 
are achieved in \cite{noi} by means of the linear contraction semigroup theory
 coupled with the abstract theory of nonhomogeneous and semilinear evolution problems; in fact,  the transmission conditions 
 (\ref{KC}) and (\ref{TC})
 allows us to prove that the linear differential operators in (\ref{sysi}) are m-dissipative and then, to apply the Hille-Yosida-Phillips Theorem (see \cite{CH}) .
The existence of global solutions when the initial data are small in 
$(H^1(\mathcal A))^2 \times H^2(\mathcal A)$ norm is proved \cite{noi} too; this result holds under  
the  further assumption  
\be\label{K} \t{ for all }\nu\in\mathcal P, \t{ for some } k\in\M^\nu, 
K_{ik}^\nu\neq 0 \t{  for all } i\in\M^\nu, i\neq k\ .
\ee

\end{section}

\bigskip

\begin{section}{\bf Non negative small stationary solutions on acyclic networks}

In this section we approach the question of existence and uniqueness of stationary solutions of problem (\ref{sysi})-(\ref{K}), with fixed mass
$$\mu:=\dsp \suM\iIi \u(x) dx\geq 0\ ,$$
in the case of an acyclic network (see Section 2).
We look for  stationary solutions $(u,v,\phi)\in (H^1(\mathcal A))^2\times H^2(\mathcal A)$. 

Obviously, the flux $v$ of a stationary solution has to be constant on each arc and has  to be null on the external arcs; in the case of acyclic networks, the boundary
and transmission conditions (\ref{vbc}), (\ref{TC}) force it to be null on each arc. In order to prove this fact we consider
an internal  arc $I_j$  and its initial  node $N_\mu$; we consider the set
$$\mathcal Q =\{ \nu\in \mathcal P: N_\nu \t{ is linked to } N_\mu  \t{ by a path not covering }  I_j\} \ 
$$
(see  Fig.1:  if, for example,  $j=9$ then    $\mu=4$, $\mathcal Q=\{1,2,3,5\}$ and the arcs in bold type  form the path which links the nodes $N_5$
and $N_4$).

\begin{figure}
\includegraphics[scale=0.3]{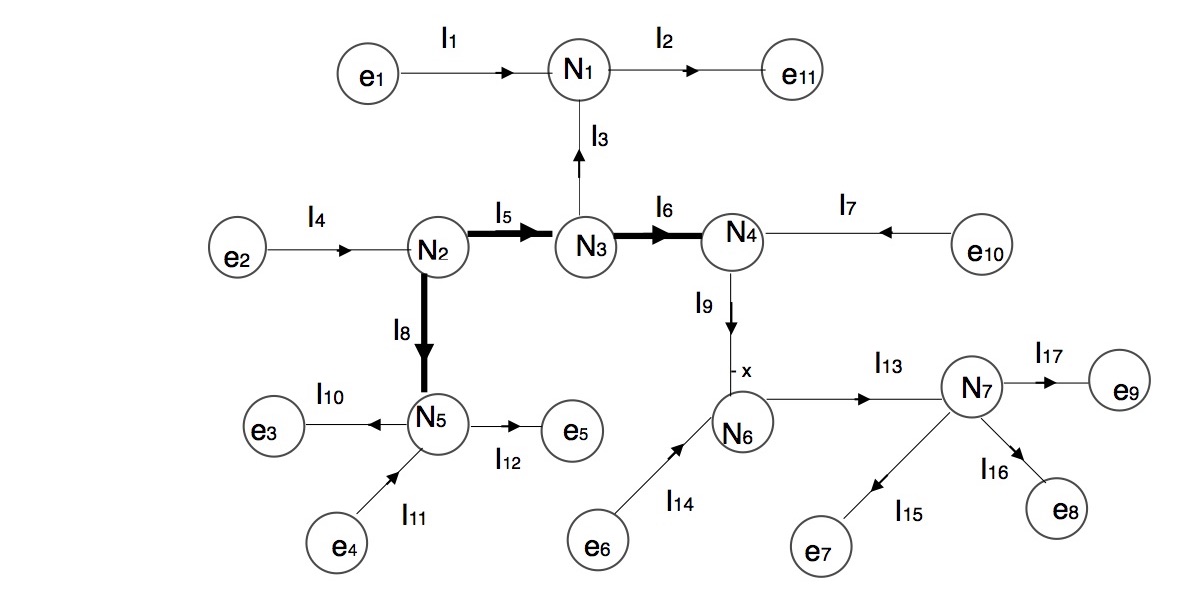}
\caption{}
\end{figure}

At each node the conservation of the flux of the density of cells, stated in Section 2,  holds; then
$$\dsp\sum_{\nu\in\mathcal Q\cup \{\mu\}} \left( \sum_{i\in I^\nu} \li\v(N_\nu) -\sum_{i\in O^\nu} \li\v(N_\nu)\right) =0\ .$$
Since, for all $i\in\M$, $v_i(x)$ is constant on  $I_i$  and $v_i(x)=0$ if $I_i$ is an external arc, the above equality reduces to
$$v_j(N_\mu)=0\ ;$$
then  $v_j(x)=0$  for all $x\in I_j$.


The previous result implies that stationary solutions must have the form $(u,0,\phi)$, where $u$ and $\phi$ have to verify the system
\be\label{stazio}\left\{\begin{array}{ll}
\li\ux=\u\phx  , \\ \\
- D_i \phxx +b_i\phii =a_i \u , 
\end{array}\right. x\in I_i,\  i\in\M,\ t>0,\ \end{equation}

with the boundary condition 
 at each outer point $e_j$, $j\in\mathcal J$,

\be\label{sbc}
{\phi_{i(j)}}_x(e_j, t)=0\qquad 
\ t>0
 \  ,\ee

and the transmission conditions, at each inner node  $N_\nu$, 

\be\label{stsu}
\suMjn K^\nu_{ij} (u_j(N_\nu)-\u(N_\nu))=0 , \qquad i\in \M^\nu\ ,
\ee

\begin{equation} \label{stcf}\begin{array}{ll}
\dsp D_i\phx(N_\nu)=\suMjn \alpha^\nu_{ij} (\phi_j(N_\nu)-\phii(N_\nu)) ,\qquad i\in \mathcal I^\nu  ,\\ \\
\dsp D_i\phx(N_\nu)=-\suMjn \alpha^\nu_{ij} (\phi_j(N_\nu)-\phii(N_\nu)) ,\qquad i\in \mathcal O^\nu .
\end{array}\ee

For each fixed inner node $N_\nu$, let $k\in\mathcal M^\nu$ be the index in condition  (\ref{K}) and let
consider the  transmission relations, for $i\in\mathcal M^\nu$,  $i\neq k$,
\be\begin{array}{ll} \dsp
0 =\sum_{j\in \M^\nu, j\neq i} K^\nu_{ij} (u_j(N_\nu)-\u(N_\nu)) \\ \\ =\dsp
\sum_{j\in \M^\nu, j\neq i,k} K^\nu_{ij} (u_j(N_\nu)-u_k(N_\nu)) -\left(\sum_{j\in \M^\nu, j\neq i} K^\nu_{ij}\right)
(u_i(N_\nu)-u_k(N_\nu))\ ;
\end{array}
\ee
the assumptions on $K^\nu_{kj}$ in (\ref{K}) ensure that the matrix of the coefficients of 
this linear system in the unknowns $(u_j(N_\nu)-u_k(N_\nu))$, $j\neq k$, is non singular (if $k=1$ it is immediate to check that it has strictly dominant diagonal).
Then the condition (\ref{stsu}) can be rewritten as
$$\dsp u_j(N_\nu)=u_k(N_\nu)\qquad \t{ for all } j\in\M^\nu.$$

Now we fix $\mu_0\geq 0$ and we look for  stationary solutions such that
\be\label{mu0}\suM \iIi\u(x) dx =\mu_0\ ;\ee
notice that for the evolution problem, the quantity $\dsp \suM \iIi\u(x,t) dx$ is preserved for all $t\geq 0$, thanks to the transmission conditions (\ref{TC}).

Integrating the first equation in (\ref{stazio}) we can rewrite   problem (\ref{stazio})-(\ref{mu0}) as the following elliptic problem on network:
\medskip

{\it 
Find $C=(C_1,C_2, ...,C_m)$
 and $\phi\in H^2(\mathcal A)$ such that}

\be\label{pr-ell2}\left\{\begin{array}{ll}\displaystyle
-D_i \phxx+b_i \phii = a_i \u \qquad x\in I_i \ ,\quad i\in\M \ ,\\ \\
\u(x)=C_i \exp\left(\frac {\phii(x)} \li\right)\qquad x\in I_i\ ,\quad i\in\mathcal M,  
\\ \\
{\phi_{i(j)}}_x(e_j) =0\ , \qquad  j\in\mathcal J,  \\ \\
\dsp D_i\phx(N_\nu)=\suMjn \alpha^\nu_{ij} (\phi_j(N_\nu)-\phii(N_\nu))\ ,\qquad i\in \mathcal I^\nu ,\quad \nu\in\mathcal P\ ,\\ \\
\dsp D_i\phx(N_\nu)=-\suMjn \alpha^\nu_{ij} (\phi_j(N_\nu)-\phii(N_\nu))\ ,\qquad i\in \mathcal O^\nu , \quad \nu\in\mathcal P\ ,\\ \\
C_j \espniu  =  C_i\espnju \ ,\qquad i,j \in\mathcal M^\nu, \quad  \nu\in\mathcal P\ ,\\ \\
\displaystyle\suM C_i\iIi \espxi dx =\mu_0\ .
\end{array}
\right.
\ee

\bigskip


We consider  the linear operator
$A:D(A)\to L^2(\mathcal A)$, 
\be \label{A2}\begin{array}{ll}
D(A)= \left\{\phi\in H^2 (\mathcal A): (\ref{sbc}), (\ref{stcf}) \t{ hold }\right\}\ ,\\ \\
A(\phi)=\left\{-D_i\phxx+b_i\phii\right\}_{i\in\M}\ ;
\end{array}\ee
then the equation in (\ref{pr-ell2}) and the boundary and transmission conditions for $\phi$ can be written as


\be A \phi = F (\phi, C) ,  \ee
where, for  $i\in\M$,  $F_i(\phi(x), C)=a_iC_i \exp\left( \frac {\phi_i(x)} {\lambda_i}\right)$.

We are going to prove the existence and uniqueness of solutions to the problem (\ref{pr-ell2}) by using the Banach Fixed Point Theorem; in order to do this we need some preliminary results about the linear equation
\be \label{0} A \phi = F (f, C^f) ,  \ee
where $f\in H^2(\A)$ is a given function, $C^f=(C^f_1,...,C^f_m)$  and $C^f_i$ are non-negative given real constants.



The existence and uniqueness of the solution $\phi\in H^2(\A)$ to the above problem (for a general $F\in L^2(\A)$ and a general network) is showed in \cite{noi}, by Lax-Milgram theorem, in the proof of Proposition 4.1; here, we need to prove some properties holding for the solution in the case of acyclic graphs, under suitable assumptions on $f$ and $C^f_i$.  

The transmission conditions (\ref{KC}) imply the following equality which will be useful in the next proofs:

\be\label{z}\begin{array}{ll}
\dsp \suM \iIi D_i \left( \phii(x)\phx(x)\right)_x \ dx= \\ \\
\dsp\sum_{\nu\in\mathcal P} \left(\sum_{i\in I^\nu} D_i \phii(N_\nu)\phx(N_\nu) -\sum_{i\in O^\nu} D_i \phii(N_\nu) )\phx(N_\nu)\right)= 
\\  \\
\dsp \sum_{\nu\in\mathcal P} \ \sum_{ij\in \M^\nu}  \alpha_{ij}^\nu  \phii(N_\nu) (\phi_j(N_\nu)-\phii(N_\nu))\ =\\ \\

 \dsp-\frac 1 2 \sum_{\nu\in\mathcal P}\  \sum_{ij\in\mathcal M^\nu} \alpha^\nu_{ij} (\phi_j(N_\nu)-\phii(N_\nu))^2 \ . 
\end{array}
\ee


Let $\dsp|\mathcal A|:=\sum_{i\in\mathcal M} | I_i|$ and $\Vert g\Vert_{\infty}:=\max \{ \Vert g_i\Vert_{\infty},  i\in\M\}$.

\bl\label{1L}
Let $\mathcal G=(\mathcal V, \A)$ be an acyclic network, let $f\in H^2(\A)$ and let $C_i^f$ be non-negative real numbers, for  $i\in \M$.
Then the solution $\phi$ to problem  (\ref{A2}),(\ref{0}) is non-negative.
Moreover, if
\be
\label{11}\suM C_i^f   \iIi\exp{ \left( \frac {f_i(x)} {\lambda_i}\right)} dx =\mu_0 \ ,
\ee
then
\be\label{phinf}\Vert \phi_{x}\Vert_{\infty}   \leq \frac {2 \max \{a_i\}_{i\in\M}}{\min \{D_i\}_{i\in\M}} 
\mu_0\\ ; \ee
 if (\ref{11}) holds and
\be\label{eta0} \dsp\Vert f_{x}\Vert_{\infty} \leq \frac {2 \max\{a_i\}_{i\in\M} 
}{\min\{D_i\}_{i\in\M}
}\mu_0 \ ,
\ee 
then there exists a quantity $K_{\mu_0}=K_{\mu_0}(a_i,b_i, D_i,\li, |\mathcal A|, \mu_0)$, depending only on the parameters  appearing in  brackets, infinitesimal when $\mu_0$ goes to zero, such that
\be \label{12}
\Vert  \phi \Vert_{W^{2,1}(\A)},  
\Vert  \phi \Vert_{H^2(\A)}
\leq  K_{\mu_0}\ .
\ee
\el

\begin{proof}

Let consider a function $\Gamma \in C^1(\R)$, strictly increasing in $(0, +\infty)$, and let $\Gamma (y)=0$ for $y\leq 0$; following
standard methods 
for   the proofs of   the maximum principle for elliptic equations  and setting $F_i(x)=F_i(f(x),C^f)$, we obtain

$$ \suM \iIi \left(-D_i (\phi_{i_x} (x)\Gamma(-\phii(x)))_x  - D_i \Gamma'(-\phii(x))\phi^2 _{i_x}(x) \, +\right.\qquad$$
$$\left. \qquad \qquad \qquad b_i\phi_i(x) \Gamma (-\phi_i(x)) -F_i(x)\Gamma (-\phii(x))\right)\ dx 
=0 \ .$$

As regard to the first term, we can argue as in (\ref{z}), taking into account the properties of $\Gamma$,

\be\label{z2}\begin{array}{ll}
\dsp \suM \iIi D_i \left(  \Gamma (-\phii)\phx\right)_x = \\ \\
 \dsp-\frac 1 2 \sum_{\nu\in\mathcal P}\  \sum_{ij\in\mathcal M^\nu} \alpha^\nu_{ij} (\phi_j(N_\nu)-\phii(N_\nu))(\Gamma (-\phi_j(N_\nu))-
 \Gamma(-\phii(N_\nu))) \geq 0 \ ;
\end{array}
\ee
the above inequality and  the non-negativity of $F_i$ imply  that 
$$\suM b_i\iIi \phii(x) \Gamma(-\phii(x))  dx \geq 0\ ,$$
so that,  thanks to the  properties of $\Gamma$ , we can conclude that $\phii(x)\geq 0$ for all $i\in\M$.



By integration of the  equation (\ref{0}), taking into account (\ref{stcf}) and (\ref{sbc}), we obtain  
\be\suM b_i\iIi  \phii(x) dx = \suM   \iIi F_i(f(x),C^f)
dx \ee
which implies
\be\label{l1} \Vert \phi\Vert_1  \leq K_1(a_i,b_i)  \mu_0\ .\ee

In order to obtain (\ref{phinf}), first we notice that, if $I_j$ is an external arc, then  the following inequality holds
$$|D_j {\phi_j}_x(x)|\leq  \int_{I_j}  D_j \vert{\phi_j}_{yy}(y)\vert \ dy
\leq \int_{I_j}\left( b_j \phi_j(y) +C_j^f a_j \exp \left(\frac{f_j(y)}{\lambda_j}\right)\right)dy \ .$$

Then we consider an  internal arc $I_j$ and   its initial node $N_\mu$ and  the sets
$$\mathcal Q =\{ \nu\in \mathcal P: N_\nu \t{ is linked to } N_\mu  \t{ by a path not covering }  I_j\} \ ,
$$
$$\mathcal S=\{ i\in \mathcal M: I_i  \t{ is incident with } N_l  \t{ for some }l \in\mathcal Q\}\ 
$$
(see Fig.1: if, for example,  $j=9$, then $\mu=4$, $\mathcal Q=\{1,2,3,5\}$, $\mathcal S=\{1,2,3,4,5,6,8,$
$10,11,12\}$);
at each node the conservation of the flux, stated in Section 2 as a consequence of the transmission conditions, holds; then
$$\dsp\sum_{\nu\in\mathcal Q\cup \{\mu\}} \left( \sum_{i\in I^\nu} D_i\phx(N_\nu) -\sum_{i\in O^\nu} D_i\phx(N_\nu)\right) =0\ .$$

Let $x$ be a point on the arc $I_j$  (see Fig.1 for  $j=9$, $\mu=4$) and $I_j^x$  the part of  $I_j$ which  connects $N_\mu$ and $x$; then, using the above equality and the boundary conditions (\ref{phibc}), we have
$$|D_j {\phi_j}_x(x)|=\left \vert D_j {\phi_j}_x(x) +\dsp\sum_{\nu\in\mathcal Q\cup \{\mu\}} \left( \sum_{i\in I^\nu} D_i\phx(N_\nu) -\sum_{i\in O^\nu} D_i\phx(N_\nu)\right) \right\vert =$$
$$\left\vert\sum_{i\in \mathcal S} \iIi D_i {\phii}_{yy}(y)\ dy + \int_{I_j^x}  D_j{\phi_j}_{yy}(y)\ dy \right\vert =
$$
$$
 \left\vert\sum_{i\in \mathcal S} \iIi \left( b_i \phii(y) -C_i^f a_i \exp \left(\frac{f_i(y)}\li\right)\right)dy +
\int_{I_j^x}\left( b_j \phi_j(y) -C_j^f a_j \exp \left(\frac{f_j(y)}{\lambda_j}\right)\right)dy \right\vert\ .$$
\medskip
Then $\dsp D_j \Vert \phi_{j_x}\Vert_{\infty}\leq 2 \mu_0 \max \{a_i\}_{i\in\M}$   for all $j\in\M$ and 
 we obtain    (\ref{phinf}) which implies

\be\label{13}
\Vert \phi_x\Vert_1 \leq 
 \frac {2 \max \{a_i\}_{i\in\M}}{\min \{D_i\}_{i\in\M}}  |\mathcal A|
\mu_0\, ,
\ee

and

\be\label{14}
\Vert \phi_{x}\Vert_2 \leq  \frac {2 \max \{a_i\}_{i\in\M}}{\min \{D_i\}_{i\in\M}}   |\A|^{\frac 1 2} \mu_0\ ;
\ee

moreover, by Sobolev embedding theorem, we obtain
$$\Vert \phi\Vert_{\infty} \leq K_4(a_i,b_i,D_i,|\mathcal A|) \mu_0\  ,$$
where $K_4$ is a suitable constant.

The estimates for the function $\phi_{xx}$ follow by using the equation (\ref{0}),
taking into account  (\ref{11}), (\ref{eta0}). 
In particular, from the equation (\ref{0}),  using (\ref{z}),
we have
$$\suM \frac {D_i^2}{b_i} \iIi\phxx^2(x) \ dx \leq
\suM \frac {\Vert F_i\Vert_{\infty}}{b_i} \iIi F_i(x)\ dx 
\leq \frac { \max \{a_i \Vert F_i\Vert_{\infty}\}}{\min\{b_i\}}\mu_0\  $$
and the embedding of $W^{1,1}(I_i) $ in $L^{\infty}(I_i)$ gives
$$ \suM \frac {D_i^2}{b_i} \iIi\phxx^2(x) \ dx \leq K_5  
(1+\Vert {f}_x\Vert_\infty)  \mu_0^2
$$
where $K_5=K_5(a_i,b_i,D_i,\li)$ .
\end{proof}

Now we  can prove the following theorem.

\bt \label{pmu2} Let $\mathcal G=(\mathcal V,\mathcal A)$ be an acyclic network. There exists $\epsilon>0$ such that, if  $ 0\leq\mu_0\leq \epsilon$ ,  then problem (\ref{sysi})-(\ref{K}) has  a unique stationary solution satisfying (\ref{mu0}); the solution has the form
$$\left(C_i \espxi, 0, \phii(x)\right)\ \ i\in \M ,$$ 
where $\phii(x)\geq 0$
and $C_i$ are nonnegative constants 
such that 
$u_j(N_\nu)=u_i(N_\nu)$  for all $\nu\in\mathcal P$, $ i,j\in \mathcal M^\nu$.\et

\begin{proof}

First we notice that, if a stationary solution $(u,v,\phi)$ satisfying (\ref{mu0}) exists, then  $u$ is non-negative, since the costants $C_i$ in (\ref{pr-ell2}) must have the same sign, so that they have to be non-negative to satisfy the condition $\mu_0 \geq 0$; arguing as in the proof of Lemma \ref{1L} we prove that $\phi$ is non-negative too. If $\mu_0>0$ then $u$ and $\phi$ are positive functions.

We are going to use a fixed point technique. Given $\phi^0 \in D(A)$, we want to define a function $u^0(x)$ on the network, such that, for $i\in\M$, 
$$u_i^0(x) =C^{\phi^0}_i \exp
\left(\frac {\phi_i^0(x)} {\lambda_i}\right)\, ,$$
 where the constants $C_i^{\phi^0}$ satisfy the following linear system composed by the last conditions in (\ref{pr-ell2})
\be\label{pr-ell3}\displaystyle
C^{\phi^0}_j \ooespniu  =  C^{\phi^0}_i\ooespnju \ ,\qquad i,j \in\mathcal M^\nu, \quad  \nu\in\mathcal P\ ,
\ee
\be\label{massa}\displaystyle\suM C^{\phi^0}_i\iIi \oespxi dx =\mu_0\ .
 \ee
 
 The system (\ref{pr-ell3}),(\ref{massa}) has a unique solution; actually, since the network has no cycles,  the system (\ref{pr-ell3}) has $\infty^1$ solutions $
 C^{\phi^0,\alpha}= (\alpha, \alpha \delta_2, \alpha \delta_3,...,\alpha \delta_m)$, $\alpha\in\R$,
 where $\delta_i$ are suitable coefficients, and the condition (\ref{massa}) determines the value of $\alpha$.
 
In order to give an explicit expression for the coefficients $\delta_i$ we consider an arc, $I_1$,
and we  define
$$u^0_1(x):=\alpha  \exp\left(\frac {\phi_1^0(x)} {\lambda_1}\right)\  .$$
Let $N_\mu$ one of the  extreme points of $I_1$, then  we define the function $u^0$ on the other arcs which are incident with $N_\mu$  in such a way to verify the equalities  in (\ref{pr-ell3}) for the node $N_\mu$,
$$u^0_j(x):=   \alpha   \exp\left(\frac {\phi_1^0(N_\mu)}{ \lambda_1}\right)\exp\left( -\frac {\phi_j^0(N_\mu)}{\lambda_j}\right)\exp\left(\frac {\phi_j^0(x)}
{\lambda_j}\right) \ \t{ for all }    j\in\mathcal M^\mu\, , j\neq 1;$$
i.e. we set $ C^{\phi^0}_j=     \alpha   \exp\left(\frac {\phi_1^0(N_\mu)}{ \lambda_1}\right)\exp\left( -\frac {\phi_j^0(N_\mu)}{\lambda_j}\right) $ , 
 $j\in\mathcal M^\mu\, , j\neq 1$.
\smallskip

\begin{figure}
\includegraphics[scale=0.3]{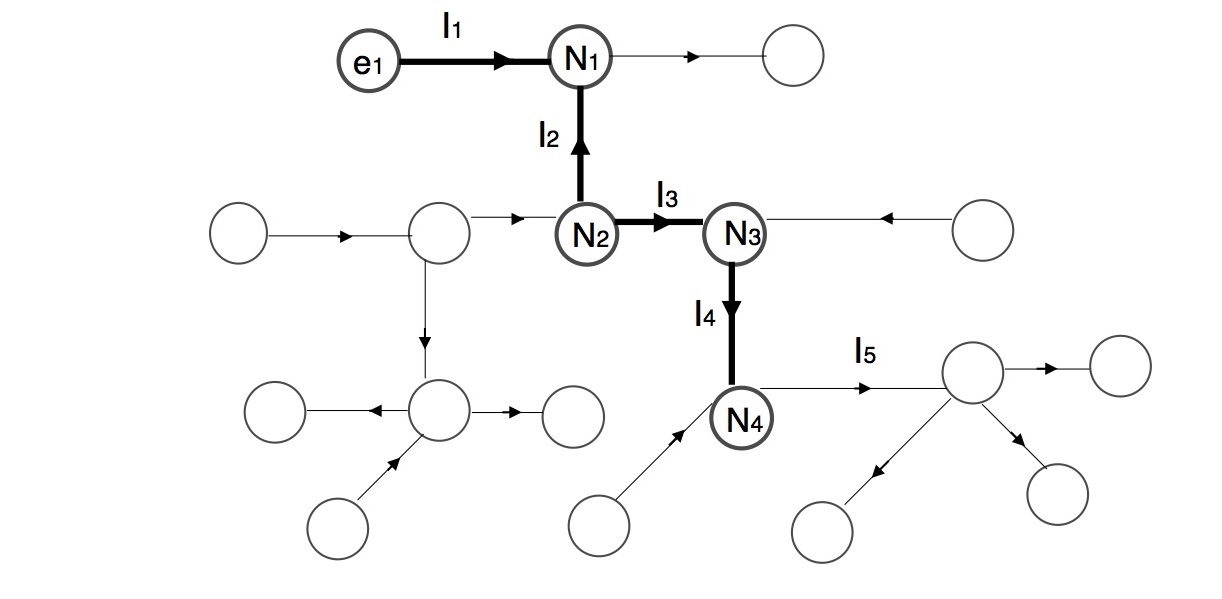}
\caption{}
\end{figure}

This procedure can be iterated at each node reached by one of the  arcs $I_j$, $j\in \M^\mu$, and at the other extreme point of $I_1$, if it is an internal arc, and so on, defining in this way  the function $u^0$ on each arc of the network. 
Notice that this construction is possible since there are no cycles in the graph.
The function $u$  can be expressed, on each arc of the network, as follows (if it is the case, 
renumbering in suitable way the arcs and the nodes):
let consider 
the path from the outer point $e_1$ to an inner node $N_{h-1}$, composed from the arcs $I_i$, $i=1,...,h-1$, (passing through the vertexes $N_i$, $i=1,...,h-1$), and let 
  $I_h$ be an arc incident with the node $N_{h-1}$, not belonging to the path (see Fig.2 where   h=5 and the highlighted arcs forms the path);
following the  procedure described before, after setting

$$\dsp\mathcal E_h(\phi^0) :=\frac {\Pi_{i=1,...,h-1} \exp \left(     \frac    {\phi^0_i(N_i)}      {\li}  \right)  }
{\Pi_{i=1,...,h-1} \exp \left(     \frac    {\phi^0_{i+1}(N_i)}      {\lambda_{i+1}}  \right)  }\ ,
$$

 we define 
$$u^0_h(x) := \alpha\mathcal E_h(\phi^0) \exp\left(\frac {
\phi^0_h(x)}      {\lambda_h}  \right)  \ .
$$

The quantity $\alpha$ is fixed in such a way to verify the last condition in (\ref{pr-ell2}),
$$\alpha\suM \mathcal E_i(\phi^0)\iIi \oespxi dx =\mu_0\ ,$$
so that, for all $i\in \M$,
\be\label{ccc} u^0_i(x)= \dsp C_i^{\phi^0} \exp\left(\frac {\phi_i^0(x)} {\lambda_i}\right)\  ,\ \ 
C_i^{\phi^0}:= \mu_0 \frac {\  \mathcal E_i(\phi^0) }
  { \dsp\  \suMj \mathcal E_j(\phi^0)  
  \int_{I_j} \oespxj dx  \ }.
\ee


Let  $G$ be the operator defined in $D(A)$  such that $\phi^1:= G(\phi^0)$ is the solution of problem
(\ref{0}) where  $f=\phi_0$ and $ C_i^f=C_i^{\phi^0}$
 for $i\in\M$,
 $$A\phi^1=F(C^{\phi^0}, \phi^0)\ ;$$
 let 
$K_{\mu_0}=K _{\mu_0}(a_i,b_i, D_i,\li,  |\mathcal A|,{\mu_0})$ be the quantity in Lemma \ref{1L}, 
and let
$$B_{\mu_0}:= \{ \phi\in D(A) : \phi\geq 0\ , 
\Vert \phi_{x}\Vert_{\infty}   \leq \frac {2 \max \{a_i\}_{i\in\M}}{\min \{D_i\}_{i\in\M}} 
\mu_0\ ,
\Vert  \phi \Vert_{H^2}
\leq K_{ \mu_0}\ \}$$
equipped with the distance $d$ generated by norm of $H^2(\A)$; $(B_{\mu_0},d)$ is a complete metric space.
From the  lemma we know that  solutions to problem (\ref{A2})-(\ref{11}) have to belong to $
B_{\mu_0}$, then $G(B_{\mu_0})\subseteq B_{\mu_0}$; next
we are proving that, if $\mu_0$ is small enough, then $G$ is a contraction  in $B_{\mu_0}$.

We consider $\phi^0,\ov\phi^0 \in B_{\mu_0}$ and the corresponding $u^0,\ov u^0$ and $\phi^1,\ov\phi^1$;
using the equation satisfied by $\phi^0$ and $\ov \phi^0$ , for all $i\in\M$ we can write
\be\label{en} \begin{array}{ll}
\displaystyle b_i \iIi (\phi^1_i(x)-\ov \phi^1_i(x))^2\, dx 
+ D_i\iIi ( \phi^1_{ix}(x)-\ov \phi^1_{ix}(x))^2 \,dx \\ \\
\displaystyle -D_i\iIi\left(( \phi^1_{ix}(x)-\ov \phi^1_{ix}(x))(\phi^1_i(x)-\ov \phi^1_i(x))\right)_x\, dx \\ \\
\displaystyle = a_i\iIi (u^0_i(x)-\ov u^0_i(x))(\phi^1_i(x)-\ov \phi^1_i(x))\ \, dx\ ;
 \end{array}\ee
 
using (\ref{z}), from (\ref{en}) we infer that 

\be\label{mai}\suM \Vert {\phii}^1 -\ov {\phii}^1\Vert_{H^2} 
\leq 
K(a_i,b_i,D_i) 
 \suM \Vert u^0_i -\ov u_i^0\Vert_{2} \ ,\ee

We set
$$ J_i^{\phi^0}:=\int_{I_i}\oespxi dx\, ,   \qquad  E_i^{\phi^0}(x):= \oespxi\ \ ;$$

we have
\be\label{ora}\begin{array}{ll} 
 \displaystyle\vert u_i^0(x)-\ov u_i^0(x)\vert =\dsp  \mu_0\left\vert \frac {\  \mathcal E_i(\phi^0) E^{\phi^0}_i(x) \ }
  { \dsp\  \suMj \mathcal E_j(\phi^0) J_j^{\phi^0} \ }
  -  \frac {\  \mathcal E_i(\ov\phi^0) E^{\ov \phi^0}_i(x) \  }{\ \dsp \suMj \mathcal E_j(\ov\phi^0) J^{\ov \phi^0}_j\ }
 \right\vert

\ .\end{array}\ee

In order to treat the above quantity we have to consider that, for all $g\in B_{\mu_0}$,
$E^g_i(x)\geq 1$, $ {J^g_i}\geq |I_i|$ and 
there exists a constant $K_6=K_6(K_{\mu_0},\li)$, increasing with  $\mu_0$, such that,
for all $i\in\M$
$$ \max_{I_i} E^g_i(x) 
\leq K_6\ ,\quad
J_i^g
\leq K_6 \vert I_i \vert\ ,$$
$$\left\vert  E^{\phi^0}_i(x) -{E^{\ov\phi^0}_i}(x)
\right\vert \leq K_6 \vert \phi_i^0(x) -\ov \phi_i^0(x) \vert\ ,$$
$$\vert J_i^{\phi^0}- J_i^{\ov \phi^0}  \vert \leq K_6 \iIi\vert \phi_i^0(x) -\ov \phi_i^0(x) \vert\, dx \ .$$

The above inequalities can be used in (\ref{ora})  so that (\ref{mai}) implies
\be\label{mai2}\suM \Vert {\phii}^1 -\ov {\phii}^1\Vert_{H^2} 
\leq \mu_0
K_7(a_i,b_i,D_i,K_{\mu_0}, |\A|) 
 \suM \Vert \phi_i^0 -\ov \phi_i^0\Vert_{H^1} \ ,\ee
 where $K_7$ increases with $\mu_0$;
hence, for $\mu_0$ small enough, $G$ is a contraction on $B_{\mu_0}$  and we can use the Banach Fixed Point Theorem.

Let $\phi$ be the unique fixed point of $G$ in $B_{\mu_0}$ and let $C^\phi=(C_1^\phi,C_2^\phi,..., C_m^\phi)$ where $C_i^{\phi}$, for $i\in\M$, are computed as in (\ref{ccc}); then $(\phi, C^{\phi})$ is the unique solution to Problem $(\ref{pr-ell2})$ and the claim is proved.

\end{proof}

For any constant $U\geq 0$, the triple $(U,0,\frac {a_i}{b_i} U)$ satisfies the equations in (\ref{sysi}) on the arc $I_i$. 
Let $Q$ be a real non-negative number; 
if $\frac {a_i}{b_i}=Q$ for all $i\in\M$, then the same triple $(U,0,Q U)$ satisfies the equations on each arc $I_i$ and it is a stationary solution to the problem (\ref{sysi})-(\ref{K}). Then, as a consequence of the previous theorem, we have the following proposition, with 
 $\ep$ as in the theorem.
 
 \bp \label{000}Let $\mathcal G=(\mathcal V,\mathcal A)$ be an acyclic network. If $\frac {a_i} {b_i}=Q$  for all $i\in \M$
  and $0\leq\mu_0\leq \epsilon$,
   then the unique stationary solution  to problem (\ref{sysi})-(\ref{K}), (\ref{mu0}) is the constant solution
$\displaystyle\left(\frac{\mu_0}{|\A|}, 0,  Q \frac{\mu_0}{|\A|}\right)$.
\end{prop}

\medskip

\noindent {\bf Remark 3.1 } 
For general networks,  when the value of $\frac {a_i}{b_i}=Q$  on each arc, the stationary solution of Proposition \ref{000} always exists.  
More precisely, if $\frac {a_i} {b_i}=Q$, in the class of the functions $(u,v,\phi)$ which are constant on each arc, the stationary solution 
$\displaystyle\left(\frac{\mu_0}{|\A|}, 0,  Q \frac{\mu_0}{|\A|}\right)$
is the unique stationary solution with  mass $\mu_0$; this fact is true without any restrictions on the value of $\mu_0$ and on the structure of the network. Actually,
if we assume that  $u$ is constant on each arc, 
then, using the equations, we infer that, on each arc,
$\phi_x(x)$ is constant too, hence
 $\phi_{xx}=0$ and $\phi(x)$ is constant. Then   $v(x)=0$ on each arc; hence, arguing as at the beginning of this section, we obtain that $u$ is continuous on the network.
 
 \medskip
 
 In the next proposition we are going to prove that, in  a set of {\it small} solutions, such stationary solution is the unique one with fixed mass $\mu_0$.
 
 \bp Let $\frac {a_i} {b_i} =Q$ for all $i\in\M$ and let $(u,v,\phi)\in (H^1(\A))^2\times H^2(\A) $  be a stationary solution of problem (\ref{sysi})-(\ref{K}),(\ref{mu0}). There exists $\ep_0>0,$ depending on $\li, a_i, b_i, D_i, \beta_i, |\A|$, such that, if $\Vert u\Vert_{H^1}\leq \ep_0$, then $(u,v,\phi)=
 \displaystyle\left(\frac{\mu_0}{|\A|}, 0,  Q \frac{\mu_0}{|\A|}\right)$. 
 \epr

\begin{proof} 
We set $H:=\Vert u\Vert_{H^1}$.
The transmission conditions (\ref{TC}) imply that 
$$ \sum_{\nu\in\mathcal P}\left( \sum_{i\in\mathcal I^\nu} \li \u(N^\nu)\v(N^\nu) - \sum_{i\in\mathcal O^\nu} \li \u(N^\nu)\v(N^\nu)\right)\geq 0\ ,$$
so, by using the first two equations in (\ref{sysi}), we obtain
$$2 \suM \beta_i \iIi v_i^2(x)dx \leq  \suM \Vert u_i\Vert_{\infty} \iIi (v_i^2(x) +\phi^2_{i_x}(x) ) \ dx \ ,$$
$$ \suM \li\iIi u^2_{i_x}(x) \ dx \leq \suM  \Vert u_i\Vert_{\infty} \iIi (u_{i_x}^2(x) +\phi^2_{i_x}(x) ) \ dx  \ +\qquad\qquad\qquad$$
$$\qquad \suM \frac {\beta_i^2}{\li}
\iIi v^2_i(x)\ dx \ ;$$
the above inequalities implies the following one
\be\label{sara} \Vert v\Vert_2^2 + \Vert u_x\Vert_2^2\leq K_0 H \left( \Vert \phi_x\Vert_2^2 + \Vert v\Vert_2^2 +\Vert u_x\Vert_2^2 \right)\ ,
\ee
where $K_0$ is a positive constant depending on the parameters $\li,\beta_i$ and the Sobolev embedding costant.

The transmission conditions (\ref{KC}) imply that 
$$ -\sum_{\nu\in\mathcal P}\left( \sum_{i\in\mathcal I^\nu} D_i\phii(N^\nu)\phx(N^\nu) - \sum_{i\in\mathcal O^\nu} D_i \phii(N^\nu)\phx(N^\nu)\right)\geq 0\ ;$$
moreover, the 
assumption (\ref{K}) imply that, for each $\nu\in\mathcal P$, for suitable coefficients $\theta^\nu_{ij}$ and suitable $k\in \M^\nu$,
$$\dsp u_j(N_\nu)=u_k(N_\nu)+\sum_{i\in \M^\nu, i\neq k} \theta^\nu_{ij} \v(N_\nu)\qquad \t{ for all } j\in\M^\nu\ ,$$
(see Lemma 5.9 in \cite{noi});  then, by the last equation in (\ref{sysi}), arguing as in the proof of Proposition 5.8 in \cite{noi}, we obtain

\be\label{sara2} \Vert \phi_x\Vert_2^2 + \Vert \phi_{xx}\Vert_2^2\leq  K_1  \left(  \Vert v\Vert_2^2 +\Vert u_x\Vert_2^2 \right)\ ,
\ee
where $K_1$ is a positive constant depending on the parameters $D_i,a_i,b_i,\theta^\nu_{ij}$.

By inequalities (\ref{sara}) and (\ref{sara2}) we deduce the following one

$$ \Vert v\Vert_2^2 + \Vert u_x\Vert_2^2\leq K_0(1+K_1) H \left(  \Vert v\Vert_2^2 +\Vert u_x\Vert_2^2 \right)\ ,
$$
which, for small $H$, implies 
$ \Vert v\Vert_2 , \Vert u_x\Vert_2=0$.

\end{proof}

In the cases when $\frac {a_i} {b_i}$ depends on the arc in consideration, stationary solutions with the component $u$ constant on each arc, can be inadmissible. As we showed before, $v$ should be zero, $u$ should be constant on the whole network and 
 $\phi$ should be constant on each arc,
$$
\u(x)=\frac {\mu_0}{|\A|} \ ,\  \ \ \phi_i(x)= \frac{a_i}{b_i} \frac{\mu_0}{|\A|}\ , \qquad i\in\M\ .$$

 Therefore the transmission conditions, for each $\nu\in\mathcal P$,
$$ \sum_{j\in\M^\nu}\alpha^\nu_{ij} \frac{\mu_0}{|\A|} \left( \frac {a_j} {b_j}-\frac {a_i} {b_i}\right)=0 \ ,\qquad i\in\M^\nu \ ,$$
are constraints on the relations between the parameters of the problem which have to hold if the  constant stationary  solution exists.

For example, in the case of two arcs, if $\frac {b_2} {a_2}\neq\frac {b_1} {a_1}$ (and $0<\mu_0\leq \ep$), the stationary  solution can not be  constant on the arcs , since the trasmission condition at the node,
$$ \alpha_{11} \frac{\mu_0}{|\A|} \left( \frac {b_2} {a_2}-\frac {b_1} {a_1}\right)=0 \ ,$$
cannot be satisfied.

Hence, in the cases when $\frac {a_i} {b_i}$ depends on the arc in consideration, if $(u,v,\phi)$ is the stationary solution  in Theorem \ref{pmu2}, then $u$ is  a continuous function on all the network but it is not constant on each arc.

\end{section}


\begin{section}{\bf Asymptotic behaviour}

In this section we are going to show that the constant stationary solutions previously introduced, provide the asymptotic profiles for  a class of solutions to problem (\ref{sysi})-(\ref{K}).  
We recall  that existence and uniqueness of global solutions 

\be\label{spazio}
\begin{array}{ll}
u,v\in C\left([0,+\infty);H^1(\mathcal A)\right)\cap C^1\left([0,+\infty);L^2(\mathcal A)\right),  \\ \\
\phi\in C\left([0,+\infty);H^2(\mathcal A)\right)\cap C^1\left([0,+\infty);L^2(\mathcal A)\right) , 

\phi_{x}\in H^1\left((0,+\infty)\times \mathcal A\right)  ,
\end{array}\ee
to such problem  is proved in \cite{noi}, when the initial data are sufficiently small in   $(H^1(\A))^2\times H^2(\A)$ norm and the following condition holds
\be\label {Q} \frac {a_i}{b_i}=Q\ \ \t{for all }i\in \M\ ;\ee 
 in particular it is proved that the functional $F$ defined by
\be\label{func}\begin{array}{ll}
\dsp F_T^2(u,v,\phi) := \suM \left(\sup_{t\in[0,T]}  \Vert \u(t)\Vert_{H^1}^2+
\sup_{t\in[0,T]}\Vert \v(t)\Vert_{H^1}^2+
\sup_{t\in[0,T]}\Vert \phx(t)\Vert_{H^1}^2\right) \\ \\
\dsp +\int_0^T 
\left( \Vert u_x(t)\Vert_{2}^2+
\Vert v(t)\Vert_{H^1}^2+\Vert v_t(t)\Vert_{2}^2+\Vert \phi_x(t)\Vert_{H^1}^2
+ \Vert \phi_{xt}(t)\Vert_{2}^2\right)\ dt \end{array}\ee
is uniformly bounded for $T>0$.

Here and below we use the notations
$$\Vert f_i(t)\Vert_2:=\Vert f_i(\cdot,t)\Vert_{L^2(I_i)} ,  \ \ 
\Vert f_i(t)\Vert_{H^s}:=\Vert f_i(\cdot,t)\Vert_{H^s(I_i)}   \ .
$$

Now we assume (\ref{Q}), we fix $\ov \mu\geq 0$ and  we consider  the constant stationary   solution,
$(\ov u,0,\ov \phi)$,
to problem 
(\ref{sysi})-(\ref{K}), 
such that
$\displaystyle \ov u |\A|= \ov \mu$; 
moreover let  
$(\tilde u_0,\tilde v_0, \tilde \phi_0)\in (H^1(\A))^2\times H^2(\A)$ be a small perturbation of  $(\ov u,0,\ov \phi)$, i.e.,
setting
$u_0:=\tilde u_0-\ov u, v_0:=\tilde v_0, \phi_0:=\tilde\phi_0-\ov \phi\ $,
the $(H^1(\A))^2\times H^2(\A)$ norm of $(u_0, v_0,  \phi_0)$ is bounded by a suitable small $\ep_0>0$.

 If $(\tilde u,\tilde v, \tilde \phi)$ is the solution to problem (\ref{sysi})-(\ref{K}) with  initial data $(\tilde u_0,\tilde v_0, \tilde \phi_0)$ and
   $u:=\tilde u-\ov u, v:=\tilde v, \phi:=\tilde\phi-\ov \phi\ $, then $(u,v,\phi)$ is solution to the system

\begin{equation}\label{ss}
    \left\{ \begin{array}{l}
\partial_t \u +\lambda_i \partial_x \v=0
\\ \\
\partial_t \v +\lambda_i \partial_x \u= (\u +\ov u)\partial_x \phi_i -\beta_i \v
\qquad\qquad x\in I_i\, ,t\geq 0\, , i\in \M,
\\ \\
\partial_t \phii= D_i \partial_{xx} \phi_i +a_i\u -b_i\phi_i\ ,
\end{array}\right.\end{equation}
complemented with the conditions (\ref{uvic})-(\ref{K}) and  initial data  $(u_0, v_0, \phi_0) $ defined above.

The existence and uniqueness of local solutions to this problem can be achieved by means of semigroup theory, following the method used in\cite{noi}, with little modifications. 

On the other hand,  if we assume that $\ov u $   is suitably small,
the method used in that paper to obtain
 the global existence result in the case of  small initial data
can be used here too, modifying the  estimates in order to treat the further term in the second equation and then using  the smallness of $\ov u$.

Below we list  a priori estimates holding for the solutions to problem (\ref{ss}), (\ref{uvic})-(\ref{compat}); we don't give the proofs since they are equal to those in \cite{noi}, in Section 5,   except for easy added calculations to treat the term $\ov u {\phi_{i}}_x$.

\bp\label{I}
Let $(u,v,\phi)$ be a local solution  to problem (\ref{ss}),(\ref{uvic})-(\ref{compat}),

$$
\begin{array}{ll}
u,v\in C\left([0,T];H^1(\mathcal A)\right)\cap C^1\left([0,T];L^2(\mathcal A)\right)\ ,  \\ \\
\phi\in C\left([0,T];H^2(\mathcal A)\right)\cap C^1\left([0,T];L^2(\mathcal A)\right)
\ ,\ 
\phi_{x}\in H^1\left((0,T)\times \mathcal A\right) \ ;
\end{array}$$
 then

\begin{itemize}
\item[a)]
$$ 
\begin{array}{ll} 
\displaystyle
\suM\left( \soT\Vert\u(t)\nld^2+\soT\Vert\v(t)\nld^2 + \beti \ioT \Vert\v(t)\nld^2 dt \right)\\ \\ \leq
C\dsp
\suM\left(\Vert\uo\nld^2+\Vert\vo\nld^2 \right)\\ \\

+\displaystyle  C\suM \left(\soT\Vert\u(t)\nhu +\ov u\right) \ioT\left(\Vert\phx(t)\nld^2 +\Vert\v(t)\nld^2
\right)dt  \ ;
\end{array}$$

\smallskip

\item[b)]
$$\begin{array}{ll}\displaystyle 
\suM\left(\soT
\Vert \vx(t)\nld^2+\soT\Vert\vt(t)\nld^2 +\ioT \Vert\vt(t)\nld^2\, dt\right)\\ \\ \leq
 \displaystyle C \left(\Vert v_0\nhu^2+\Vert u_0\nhu^2\Vert \phi_0\nhd^2\right)
\\ \\ +
 \displaystyle C \suM
\left(\soT\Vert \u(t)\nhu + \ov u \right)\ioT \left({\phii}_{xt}(t)\nld^2 +\Vert {\v}_t(t)\nld^2\right)\ dt \\ \\
+
 \displaystyle C \suM\soT\Vert \phi_{x}(t)\nhu
\ioT\left( \Vert
{\v}_t(t)\nld^2+\Vert\v(t)\nhu^2\right)\, dt \ ;
\end{array}\ $$

\smallskip

\item[c)]
$$\begin{array}{ll}\displaystyle 
\suM\soT\Vert \ux(t)\nld^2\leq C\suM\left(\soT\Vert\vt(t)\nld^2+\soT\Vert\v(t)\nld^2\right)\\ \\
\displaystyle
+C \suM\left(\soT\Vert \u(t)\nhu+\ov u\right)\left(\soT\Vert\ux(t)\nld^2 +\soT\Vert \phx(t)\nld^2\right)
\ ;
\end{array}$$

\smallskip

\item[d)]
$$\begin{array}{ll}\displaystyle
\suM\ioT\Vert \ux(t)\nld^2\, dt\leq C\suM\ioT\left(\Vert\vt(t)\nld^2+\Vert\v(t)\nld^2\right)
\,dt
\\ \\
\displaystyle
+C \suM\left(\soT\Vert \u(t)\nhu +\ov u\right)
\ioT\left(\Vert\ux(t)\nld^2 +\Vert \phx(t)\nld^2\right)\,dt 
\ ;
\end{array}$$

\item[e)]
$$\begin{array}{ll}\displaystyle
 \suM\ioT\Vert \vx(t)\nld^2 \,dt\leq 
C \suM\left(\Vert \vo\nld^2+ \Vert \uo\nhu^2\left(1+\Vert \pho\nhu^2\right)
\right)\\ \\ \dsp  
+C \suM\left(\ioT\Vert v_{it}(t)\nld^2\,dt+\soT\Vert \vt(t)\Vert_2^2\right)
\\ \\
\displaystyle
+C\suM 
\left(\soT\Vert \u(t)\nhu+\soT\Vert\phx(t)\nhu +\ov u\right) \\ \\
\qquad\qquad\qquad \displaystyle
\times\ioT\left(\Vert \v(t)\nhu^2 + \Vert \phi_{ixt}(t)
\nld^2\right)\,dt\ ;
\end{array}$$

\smallskip

\item[f)]
$$\begin{array}{ll}\dsp
\suM\left( \soT \Vert \pht(t)\nld^2 +\ioT\left(\Vert \pht(t)\nld^2+\Vert \phtx(t)\nld^2\right) \ dt\right) \\ \\ \qquad\qquad\leq 
\dsp
C\suM \left(\Vert \pho\nhd^2+\Vert \uo\nld^2+\ioT \Vert \ut(t)\nld^2\right)\ 
\end{array}$$

\smallskip

\item[g)]
$$\begin{array}{ll}\dsp
\suM\left( \soT\Vert \phxx(t)\nld^2 +\soT\Vert \phx(t)\nld^2\right)\\ \\ \qquad\qquad \leq 
\displaystyle C\suM \left( \soT\Vert \pht(t)\nld^2 +\soT\Vert \u(t)\nld^2\right) \ ;
\end{array}\ $$

\smallskip
 
\item[h)]
  if (\ref{K}) and (\ref{Q}) hold,
  then

$$\begin{array}{ll}\displaystyle
\suM\ioT\left( \Vert \phx(t)\nld^2+\Vert \phxx(t)\nld^2\right)\ dt \\ \\\qquad\qquad \leq 
\displaystyle
C \suM
\ioT
\left(  \Vert\ux(t)\nld^2 +
  \Vert \v(t)\nhu^2 +\Vert {\phii}_{t}(t)\nld^2
   \right)\  dt\  ,
\end{array}$$

\end{itemize}

for suitable costants $C$.
\epr

\medskip

The estimates in the previous proposition allow to prove the following theorem about the existence of global solutions to problem  
(\ref{ss}),(\ref{uvic})-(\ref{K}),

Let $F_T(u,v,\phi)$ be the functional defined in (\ref{func}) .

\bt\label{ge} 
Let (\ref{Q}) hold.
 There exists $\ep_0,\ep_1>0$ such that, if 
$$\ov u\leq \ep_1\ ,\ 
\ \Vert u_0\Vert_{H^1},\Vert v_0\Vert_{H^1}, \Vert \phi_{0}\Vert_{H^2}\leq \ep_0 \ ,$$
then there exists a unique global solution $(u,v,\phi)$ to problem (\ref{ss}),(\ref{uvic})-(\ref{K}), 
$$u,v \in C([0,+\infty);H^1(\A))\cap C^1([0,+\infty);L^2(\A)) \ , $$
$$\phi\in C([0,+\infty);H^2(\A))\cap C^1([0,+\infty);L^2(\A))\ , \ 
\phi_{x}\in H^1\left( \mathcal A\times (0,+\infty)\right) 
\ .$$
Moreover, $F_T(u,v,\phi)$ is bounded, uniformly in $T$.
\et

\begin{proof}

It is sufficient to prove that the functional $F_T(u,v,\phi)$  is bounded, uniformly in $T$.

We notice that  each term in $F_T^2(u,v,\phi)$ is in the left  hand side of one of the estimates in Proposition \ref{I}, therefore, arranging all the estimates,
 we can prove the following inequality 
$$F_T^2(u,v,\phi)\leq c_1 F_0^2(u,v,\phi) +c_2 \ov u F^2_T(u,v,\phi) +c_3 F_T^3(u,v,\phi)\ ,$$
taking into account  also that, on the right hand side of the estimates, the quadratic terms (not involving  initial data) which have not the coefficient
$\ov u$, can be bounded by means of cubic ones.

If $\ov u$ is sufficiently small, the previous inequality implies
$$F_T^2(u,v,\phi)\leq c_4 F_0^2(u,v,\phi)  +c_5 F_T^3(u,v,\phi)$$
for suitable positive constants $c_i$.

It is easy to verify that, if $y_0$ is a  sufficiently small positive real number and $ h(y)=c_5 y^3-y^2 +c_4 y_0$
then there exists $0<\ov y<\frac 2 {3 c_5}$ such that $h(y)>0$ in $[0,\ov y)$ and $h(y)<0$ in $(\ov y,\frac 2 {3 c_5}]$ .

Then we can conclude that, if $F_0$ is suitably small , then $F_T$ remains bounded for all $T>0$.

\end{proof}

\bigskip

The above result, in particular the uniform, in time,  boundedness of the functional $F_T$, allow us to prove the  theorem below. 

Let 
(\ref{Q}) hold and let $(\ov u,0,\ov \phi)$ be the constant stationary solution to problem (\ref{sysi})-(\ref{K}) such that  $ \ov u |\A|=\ov \mu
$; moreover, let $\mathcal C(\A)$ be the set of  the funcions $f$ defined on $\A$ such that $f_i\in\mathcal C(I_i)$ for $i\in \M$.

\bt\label{ab}
Let (\ref{Q}) hold.
There exist $\ep_0, \ep_2>0$ such that, if  $\ov u \leq \ep_2$, 
$\displaystyle\suM\iIi  u_0(x)= \ov \mu$
and
$
\Vert( u_0-\ov u,  v_0,  \phi_0-\ov \phi)\Vert_{(H^1)^2\times  H^2} \leq \ep_0$,
then
problem (\ref{sysi})-(\ref{K}) has  a unique global solution $(u,v,\phi)$ ,
$$u,v \in C([0,+\infty);H^1(\A))\cap C^1([0,+\infty);L^2(\A)) \ , $$
$$\phi\in C([0,+\infty);H^2(\A))\cap C^1([0,+\infty);L^2(\A))\ ,$$
and, for all $i\in\M$,
$$\dsp\lim_{t\to+\infty}\Vert \u(\cdot,t) -\ov u\Vert_{C(\ov I_i)}\ ,
\lim_{t\to+\infty}\Vert  \v(\cdot,t)\Vert_{C(\ov I_i)}\ ,
\lim_{t\to+\infty}\Vert  \phii(\cdot,t) -\ov \phi\Vert_{C^1(\ov I_i)} =0\ .
$$
\et

\begin{proof}

Let $(u,v,\phi)$ be the local solution to problem (\ref{sysi})-(\ref{K}) and 
let 
$$\hu:=u-\ov u, \hv:=v,\hf:=\phi-\ov\phi ;$$
 we already noticed that 
$(\hat u,\hat v,\hat \phi)$  is  the local solution to system (\ref{ss}) complemented   by the  initial condition $(u_0-\ov u, v_0,\phi_0-\ov \phi)$ and the same boundary and  transmission condition given by (\ref{vbc})-(\ref{K}) for system (\ref{sysi}).

For suitable $\ep_0,\ep_2$ the  assumptions of Theorem \ref{ge}  are satisfied, then we obtain   the uniform boundedness of the functional 
$F_T(\hu,\hv,\hf)$, for $T\in [0,+\infty)$.

Hence the set $\{\hu(t),\hv(t),\hf(t)\}_{t\in[0,+\infty)}$ is uniformly bounded in $(H^1(\A))^2\times H^2(\A)$; thus,  if we call $E_s$ the set of accumulation points of $\{\hu(t),\hv(t),\hf(t)\}_{t\geq s}$ in $(C(\A))^2\times C^1(\A)$, then $E_s$ is not empty and  $\dsp E:=\cap_{s\geq 0} E_s \neq \emptyset$.

Let $(U(x), V(x),\Phi(x))$  be such that, for a sequence $t_n\to +\infty$, 

\be\label{convergence}\begin{array}{ll} \dsp\lim_{n\to +\infty} \suM\Vert \hu_i(\cdot, t_n)- U_i(\cdot)\Vert_{C(\ov I_i)} =0\ ,\\ \\
\dsp\lim_{n\to +\infty} \suM\Vert \hv_i(\cdot, t_n)- V_i(\cdot)\Vert_{C(\ov I_i)}  =0\ ,\\ \\
\dsp\lim_{n\to +\infty}\suM \Vert \hf_i(\cdot, t_n)- \Phi_i(\cdot)\Vert_{C^1(\ov I_i)}= 0\ .
\end{array}\end{equation}

In order to identify the limit functions we notice that $\dsp\suM \iIi U_i(x) dx =0$,
since $\dsp \suM \iIi \hu(x,t_n)\, dx=0$ for all $t_n$.

Moreover, since $\dsp \hv_i\in H^1(I_i\times (0,+\infty))$  for all $i\in\M$, if we set 
$$\omega_i(t):=\Vert \hat v_i(t,\cdot)\Vert_{L^2(I_i)}$$ then 
$\omega_i\in H^1((0,+\infty))$  and, as a consequence, $\dsp \lim_{t\to+\infty} \omega_i(t) =0$.

As 
$\dsp\lim_{n\to +\infty}\Vert \hv_i(\cdot,t_n)\Vert_2=\Vert V _i(\cdot)\Vert_2$, we obtain $\Vert V\Vert_2 = 0$.

The same argument can be applied to the functions $\hf_{i_x}$ since they belongs to $H^1(I_i\times (0,+\infty))$.  Finally, it can be applied  to  the functions $a_i\hu_i -b_i\hf_i$ since $\hat \phi_{i_t},\hat  \phi_{i_{xx}}, \hat u_{i_x}, \hat \phi_{i_x}\in L^2(I_i\times (0,+\infty))$, thanks to the uniform boundedness of $F_T(\hu,\hv,\hf) $
and to estimate f) in Proposition  \ref{I}.

As a consequence we have that 
$$ V_i(x)=0,\ \quad 
 a_i U_i(x)- b_i   \Phi_i(x)=0\ ,\ \qquad
\Phi_i(x)= \ov \Phi_i, \qquad\ x\in I_i\ ,$$
where $\ov \Phi_i$ are real numbers, so that the limit function is given by $(\frac {b_i}{a_i} \Phi_i,0, \Phi_i)$ in each interval $I_i$, for all $i\in\M$. It is easily seen that such function is a stationary solution  to problem (\ref{sysi})-(\ref{K}), which is  constant in each arc $I_i$; in particular it verifies the transmission conditions  since $(\hu,\hv,\hf)$ verifies them and the convergence result (\ref{convergence}) holds.

The condition $\dsp\suM \iIi U_i(x) dx =0$ and Remark 3.1 imply that $\Phi_i=0$ for all $i\in\M$, so that 
we can conclude that  the unique possible limit function is $(U(x),V(x),\Phi(x))=(0,0,0)$; this fact proves the claimed convergence results .

\end{proof}


\end{section}

\end{document}